\documentclass[review]{elsarticle}

\usepackage{lineno,hyperref}
\modulolinenumbers[5]

\journal{Journal of \LaTeX\ Templates}

\usepackage{amsfonts}
\usepackage{amsmath,amsthm,amssymb,mathabx}
\usepackage{booktabs}
\usepackage{xcolor}
\usepackage{color}
\usepackage{bm}
\usepackage{IEEEtrantools}
\usepackage{txfonts,pxfonts}
\usepackage[ruled]{algorithm2e}
\usepackage{algpseudocode}
\usepackage{multirow}
\usepackage{rotating}
\usepackage{soul}
\usepackage{float} % para usar [H]
\usepackage{rotating}
\usepackage{eurosym}
\usepackage{array}
\newcolumntype{L}[1]{>{\raggedright\let\newline\\\arraybackslash\hspace{0pt}}m{#1}}
\newcolumntype{C}[1]{>{\centering\let\newline\\\arraybackslash\hspace{0pt}}m{#1}}
\newcolumntype{R}[1]{>{\raggedleft\let\newline\\\arraybackslash\hspace{0pt}}m{#1}}

\def\R{\mathbb{R}}

\def\Kc{\mathcal{K}_C}
\def\Kcp{\mathcal{K}_{C}^{+}}

\newtheorem{definition}{Definition}
\newtheorem{theorem}{Theorem}
\newtheorem{corollary}{Corollary}

\newtheorem{proposition}{Proposition}

\newtheorem{remark}{Remark}

\newcommand{\splitatcommas}[1]{%
  \begingroup
  \ifnum\mathcode`,="8000
  \else
    \begingroup\lccode`~=`, \lowercase{\endgroup
      \edef~{\mathchar\the\mathcode`, \penalty0 \noexpand\hspace{0pt plus 1em}}%
    }\mathcode`,="8000
  \fi
  #1%
  \endgroup
}

\newcommand{\tuple}[1]{(\splitatcommas{#1})}

\begin{document}
%%%%%%%%%%%%%%%%%%%%%%%

\begin{frontmatter}

\title{A novel slacks-based interval DEA model \\ and application.
%to sustainable tourism in the Mediterranean region.
%\tnoteref{mytitlenote}}
%\tnotetext[mytitlenote]{Fully documented templates are available in the elsarticle package on \href{http://www.ctan.org/tex-archive/macros/latex/contrib/elsarticle}{CTAN}.
}

\author[cadiz]{Manuel Arana-Jim\'enez %\fnref{Arana}
}
\ead{manuel.arana@uca.es}
\author[cadiz]{Julio Lozano-Ram\'{\i}rez
}
%\fntex[Sanchez]{Corresponding author}
\ead{julio.lozano@uca.es}

%\fntex[Sanchez]{Corresponding author}
\author[cadiz]{M. Carmen S\'anchez-Gil \fnref{Sanchez}}
\ead{mcarmen.sanchez@uca.es}
\address[cadiz]{Department of Statistics and Operational Research, University of C\'adiz, Spain}

\author[cadiz]{Atefeh Younesi %\fnref{Younesi}
}

\ead{atefeh.younesi@uca.es}

\author[sevilla]{Sebasti\'an Lozano}
\address[sevilla]{Department of Industrial Management, 
University of Seville, Spain}
\ead{slozano@us.es}

%\author[mymainaddress,mysecondaryaddress]{Department of Statistics and Operational Research}

\cortext[mycorrespondingauthor]{mcarmen.sanchez@uca.es}

%%%%%%%%%%%%%%%%%%%%%%%
\begin{abstract}
This paper proposes a novel slacks-based interval DEA approach that computes interval targets, slacks, and crisp inefficiency scores. It uses interval arithmetic and requires solving a mixed-integer linear program. The corresponding super-efficiency formulation to discriminate among the efficient units is also presented. 
We also provide a case study of its application to sustainable tourism in the Mediterranean region, assessing the sustainable tourism efficiency of twelve Mediterranean regions to validate the proposed approach. The inputs and outputs cover the three sustainability dimensions and include GHG emissions as an undesirable output. Three regions were found inefficient, and the corresponding inputs and output improvements were computed. A total rank of the regions was also obtained using the super-efficiency model.
\end{abstract}

\begin{keyword}
\texttt{Efficiency; interval data; slacks-based model; super-SBI; GHG emissions}
%\sep \LaTeX\sep Elsevier \sep template
\MSC[2018] %90C70\sep 03E72\sep 90C29
\end{keyword}
%%%%%%%%%%%%%%%%%%%%%%%
\end{frontmatter}

%\linenumbers

\section{Introduction}\label{sec1}

Data Envelopment Analysis (DEA) is a well-known non-parametric methodology to evaluate the efficiency of a set of Decision-Making Units (DMUs) that consume inputs (i.e., resources) to produce outputs (Zhu \cite{Zhu02}, Cooper et al. \cite{Cooper2004}). A Production Possibility Set (PPS) is derived from the observed inputs and outputs using certain axioms and the Principle of Minimum Extrapolation. This PPS contains all the feasible operating points. The non-dominated subset of the PPS is called the efficient frontier. A DMU is efficient if it lies on the efficient frontier. Otherwise, it is inefficient and can be projected onto a target operating point on the efficient frontier. 

Regarding interval DEA, there exists extensive literature, most of them containing radial multiplier formulations (e.g., Despotis and Smirlis \cite{despotis}, Zhu \cite{zhu}). Also, there exist additive imprecise DEA approaches (e.g., Lee et al. \cite{lee}), FDH interval DEA models (e.g., Jahanshaloo et al. \cite{Jahanshahloo04}), non-radial, non-oriented imprecise DEA approaches (e.g., Azizi et al. \cite{azizi}), Ideal point approaches (e.g., Jahanshahloo et al. \cite{jahanshahloo1}), inverted DEA approaches (e.g., Inuiguchi and Mizoshita \cite{Inuiguchi}), interval DEA with negative data (e.g., Hatami-Marbini et al. \cite{hatami3}), flexible measure interval DEA approaches (e.g., Kordrostami and Jahani Sayyad Noveiri \cite{kordrostami3}), common weights imprecise DEA approaches (e.g., Hatami-Marbini et al. \cite{hatami2}), three-Stage DEA Model with Interval Inputs and outputs (e.g., Cheng et al. \cite{cheng}), Two stage interval DEA (e.g., Kremantzis et al. \cite{kremantzis}. Applications include the manufacturing industry (e.g., Wang et al. \cite{Wang05}), banks and bank branches (e.g., Jahanshaloo et al. \cite{jahanshahloo2}, Inuiguchi and Mizoshita \cite{Inuiguchi}, Hatami-Marbini et al. \cite{hatami3}), power plants (e.g., Khalili-Damghani et al. \cite{khalili}), etc.

Regarding super-efficiency, aimed at discriminating between efficient DMUs, it was introduced by Andersen and Petersen \cite{Andersen01} and has been continued by many authors such as Zhu \cite{Zhu01}, Seiford and Zhu \cite{seiford},  Zhong et al. \cite{zhong}, Li et al. \cite{li}, Esteve et al. \cite{esteve}, Bolos et al. \cite{bolos}, among others. In comparison, the initial super-efficiency approaches involved radial DEA models, a super-slacks-based measure of efficiency (super-SBM, \cite{tone}), and a super Slack-based measure of inefficiency (super-SBI, \cite{Moreno01}) have also been proposed.

This paper considers inputs and outputs with continuous interval-valued data as mathematical uncertainty modeling. The most similar and recent DEA model is that by Arana-Jimenez et al. \cite{arana2021}. In that work, authors consider a hybrid scenario in which, in addition to integer interval data, some inputs or outputs are given as continuous intervals, but no super-efficiency model is proposed to rank or discriminate the efficient DMUs. In the present work, we focus on continuous cases, with an enhanced formulation of the model given in \cite{arana2021} and a super-efficiency model. Thus, while \cite{arana2021} uses two phases for the slacks-based model, under an additive and non-oriented approach, we propose a one-phase interval slack-based model in the present work. This is possible by using the $gh$-difference of two intervals. Also, we propose a super SBI interval model to discriminate between efficient units. 

To illustrate the usefulness of the proposed approach, we present an application to sustainable tourism. The first definitions of sustainable tourism focused on environmental and economic development, with community involvement included later \cite{Hardy}. The current sustainable tourism policy is often economic-growth oriented, with theoretical differences from sustainable development \cite{Guo et al}. However, sustainable tourism management policies should maximize economic benefits while minimizing adverse environmental impacts \cite{Nepal}. To the latter, GHG emissions are considered undesirable output, discussed in the application section and included in our model. To this matter, Tone, and Tsutsui \cite{tone09} propose a model in which undesirable output variables are treated as inputs, as well as the cases of \cite{Ma01} and \cite{Chen01}, considering water consumption and tourism energy, respectively. Finally, tourism sustainability should not only be based on a guide of good practices; it must also consider quantitative data that allows evaluation with which to make decisions. In this regard, using indicators is a crucial tool that enables an assessment of the transition toward Sustainability \cite{Alfaro1}.

The structure of the paper is the following. Section 2 introduces the conventional crisp production possibility set (PPS) and a slacks-based DEA model for efficiency assessment. Section 3 introduces continuous intervals, especially arithmetic operations, and partial orders. Section 4 presents the continuous interval PPS and Slack-based measures of inefficiency. In Section 5, a new Enhanced interval slacks-based DEA approach is proposed. Section 6 presents the corresponding Super SBI interval DEA model to discriminate among efficient units. Section 7 presents the application to sustainable tourism, and finally, Section 8 summarizes and  concludes.

\section{Crisp production possibility set and Slack-based measure}

Let us consider a set of $n$ DMUs. For $j\in J=\{1,\ldots,n\}$, each $DMU_j$  has $m$ inputs 
$X_j = (x_{1j},\ldots,x_{mj}) \in \R^m$,  
produces $s$ outputs 
$Y_j = (y_{1j},\ldots,y_{sj}) \in \R^s$. In the classic Charnes et al. \cite{charnes78} DEA model, the production possibility set (PPS) or technology, denoted by $T$, satisfies the following axioms:
\begin{itemize}
    \item [(A1)] Envelopment: $(X_j,Y_j)\in T$, for all $j\in J$.
     \item [(A2)] Free disposability: $(x,y)\in T$, $(x',y')\in\R^{m+s}$, $x'\geqq x$, $y'\leqq y\Rightarrow (x',y')\in T$.
    \item [(A3)] Convexity: $(x,y), (x',y')\in T$, then $\lambda (x,y)+(1-\lambda)(x',y')\in T$, for all $\lambda\in[0,1]$.
    \item [(A4)] Scalability: $(x,y)\in T\Rightarrow (\lambda x,\lambda y)\in T$, for all $\lambda\in\R_+$.
\end{itemize}
%\\
Following the minimum extrapolation principle (see \cite{banker84}), the DEA PPS, which contains all the feasible input-output bundles, is the intersection of all the sets that satisfy axioms (A1)-(A4) or axioms (A1)-(A3) can be expressed, respectively, as
$$T_{DEA} = \left\{(x,y)\in  \R_{+}^{m+s} : x \geq \sum_{j=1}^{n}\lambda_j X_j, 
y \leq \sum_{j=1}^{n}\lambda_j Y_j, \lambda_j\geq 0 \right\}$$
$$T_{DEA} = \left\{(x,y)\in  \R_{+}^{m+s} : x \geq \sum_{j=1}^{n}\lambda_j X_j, 
y \leq \sum_{j=1}^{n}\lambda_j Y_j, \sum_{j=1}^{N}\lambda_j=1, \lambda_j\geq 0 \right\}.$$
%\\
In the first case, the DEA technology considers Constant Returns to Scale (CRS), while in the second case, the DEA technology considers Variable Returns to Scale (VRS).
Let us also recall that a DMU $p$ is said to be efficient if and only if for any $({x},{y})\in T_{DEA}$ such that ${x}\preceqq {X}_p$ and ${y}  \succeqq {Y}_p$, then $({x},{y})=({X}_p,{Y}_p)$. This can be determined by solving the following normalized slacks-based DEA model
\begin{eqnarray}\label{prob: CIDEA}
\mbox{(DEA)}\ \ {I}(X_p, Y_p)= 
&\mbox{Max} & {\displaystyle\sum_{i=1}^{M} \frac{s_i^x}{x_{ip}} +
\sum_{r=1}^{S} \frac{s_r^y}{y_{rp}}}
\\[0.25em]
&\mbox{s.t.} & {\displaystyle\sum_{j=1}^{N}} \lambda_j {x}_{ij} \leq {x}_{ip} - s_i^x ,\quad i=1,\dots, M,
\nonumber\\[0.25em]
&& {\displaystyle\sum_{j=1}^{N}} \lambda_j {y}_{rj} \geq {y}_{rp} + s_r^y ,\quad r=1,\dots, S,
\nonumber\\[0.25em]
&& \lambda_j \geq 0 ,\quad j=1,\dots, N, \nonumber \\[0.25em]
&& ({\displaystyle\sum_{j=1}^{N}} \lambda_j=1)\nonumber \\[0.25em]
&& s_i^x ,s_r^y \geq 0 , \quad i=1,\dots, M, \ \ r=1,\dots, S. \nonumber
\end{eqnarray}
where $\lambda_j$, $j=1,\dots,n$, are the intensity variables used for defining the corresponding efficient target of $DMU_p$ and the constraint $\sum_{j=1}^{N}\lambda_j=1$ only applies to the VRS case, not to the CRS case.
The inefficiency measures $I(X_p, Y_p)$ is units invariant and non-negative. Moreover, a $DMU_p$ is efficient if and only if $I(X_p, Y_p)=0$.

\section{Notation and preliminaries}\label{sec2}
%%%%%%%%%%%%%%%%%%%%%%%
This paper presents uncertainty on the production possibility set by modeling the corresponding inequality relationships using partial orders on integer intervals. 
This requires introducing first the following notation and results.

Let $\R$ be the real number set. We denote by $\mathcal{K}_{C}=\left\{ \left[ \underline{a},\overline{a}\right] \;|\;%
\underline{a},\overline{a}\in \mathbb{R}\mbox{ and
}\underline{a}\leq \overline{a}\right\}$ the family of all bounded closed intervals in $\R$. Some useful and necessary arithmetic operations for the purpose of this manuscript are described following (see, for instance, \cite{Stefanini2009,Stefanini2010}).
%%%%%%%
\begin{definition}\label{def interval arithmetics}
Let $A=[\underline{a},\overline{a}] \in \mathcal{K}_{C}$, and $B=[\underline{b},\overline{b}] \in \mathcal{K}_{C}$ 
\begin{itemize}
\item  Addition: 
$A+B:= \{a+b\mid a \in A , b \in B \}=[\underline{a}+\underline{b},\overline{a}+\overline{b}],$
\item  Opposite value: 
$-A=\{-a:a \in A\}=[-\overline{a},-\underline{a}],$
\item  Substraction: 
$A-B:= \{a-b\mid a \in A , b \in B \}=[\underline{a}-\overline{b},\overline{a}-\underline{b}],$%
\item  Multiplication:
$A \cdot B:= \{a \cdot b\mid a \in A , b \in B \}=
[min\{ \underline{a} \cdot \underline{b}, \underline{a} \cdot \overline{b}, \overline{a} \cdot \underline{b},\overline{a} \cdot \overline{b} \} \ , \ max\{ \underline{a} \cdot \underline{b}, \underline{a} \cdot \overline{b}, \overline{a} \cdot \underline{b},\overline{a} \cdot \overline{b} \} ]
$.
\item  Multiplication by a scalar: for any $\lambda \geq 0$, 
$\lambda \cdot A = [\lambda \cdot \underline{a}, \lambda \cdot \overline{a}] $.  
For $\lambda < 0$, 
$\lambda \cdot A = [\lambda \cdot \overline{a}, \lambda \cdot \underline{a}]$.
\end{itemize}
\end{definition}
Note that $A-A\ne 0$, in general. To overcome this issue, when the modeling requires it, we have the \emph{$gH$-difference} of two intervals $A$ and $B$, which we recall from 
\cite{Stefanini2009,Stefanini2010}, as follows:
\begin{equation}
\label{eq:gH diff}
A\ominus_{gH}B=C \Longleftrightarrow
\left\{
\begin{array}{ll}
&(a)\ A= B+C,\\
or&(b)\ B = A + (-1)C.
\end{array}
\right.
\end{equation}
Note that the difference between an interval and itself is zero, that is,  $A\ominus_{gH}A=[0,0]$.
Furthermore, the \emph{$gH$-difference} of two intervals always exists and is equal to
\begin{equation}
\label{eq:gH diff bis}
A\ominus_{gH}B = [\min\{\underline{a}-\underline{b}, \overline{a}-\overline{b}\}, \max\{\underline{a}-\underline{b}, \overline{a}-\overline{b}\}]\subset A-B.
\end{equation}
Later, in the next section, we discuss the modeling with slack variables, and we will recover the difference between intervals to that matter.
It is also necessary to define a partial order relationship for integer intervals. To this aim, we will adapt LU-fuzzy partial orders on intervals, which are well-known in the literature,
(see, e.g., \cite{wu2009b, Stefanini2009} and the references therein).
%%%%%%%%%%%%%%%%%%%%%%%%%%%
\begin{definition}\label{def interval partial orders}
Given two intervals  $A=[\underline{a},\overline{a}], B=[\underline{b},\overline{b}]\in\Kc$, we say that:
	\begin{itemize}
\item[(i)] $[\underline{a},\overline{a}] \preceqq [\underline{b},\overline{b}]$ if and only if $\underline{a} \leq \underline{b}$ and $\overline{a} \leq \overline{b}$.
\item[(ii)] $[\underline{a},\overline{a}] \prec [\underline{b},\overline{b}]$ if and only if $\underline{a} < \underline{b}$ and $\overline{a} < \overline{b}$.
\item[(iii)] $[\underline{a},\overline{a}] \precneqq[\underline{b},\overline{b}]$ if and only if $\underline{a} < \underline{b}$ and $\overline{a} \leq \overline{b}$, or $\underline{a} \leq \underline{b}$ and $\overline{a} < \overline{b}$.
	\end{itemize}
\end{definition}
Similarly, we define the relationships $A \succeqq B$ and $A\succ B$ for intervals, meaning  $B\preceqq A$ and $B\prec A$, respectively. We denote $\Kcp$ as the set of all non-negative intervals, that is, $\Kcp=\{A\in \Kc: A \succeqq 0 \}$.

%%%%%%%%%%%%%%%%%%%%%%%%%%%%%
\section{Proposed interval PPS and Slack-based measure of inefficiency}\label{sec3}

Let us consider a set of $N$ DMUs, $DMU_j$ for $j\in \{1,\ldots,N\}$, with  $M$ inputs $X \in (\Kcp)^{M}$, with $x_{ij} = [\underline{x_{ij}},\overline{x_{ij}}]\in\Kcp$ for all $i \in O^X=\{1,\dots, M\}$, and $S$ outputs $Y \in (\Kcp)^{S}$, with $y_{rj}=[\underline{y_{rj}},\overline{y_{rj}}]\in\Kcp$ for all $r\in O^{Y}=\{1,\dots, S\}$.  
And the following axioms, are analogous to (A1)-(A4) in Section 2 but consider interval inputs and outputs. We apply the corresponding interval arithmetic, Definition \ref{def interval arithmetics}, and partial order introduced in Definition \ref{def interval partial orders}:

\begin{itemize}
    \item [(B1)] \textit{Envelopment: $({X}_j,{Y}_j)\in T$, for all $j\in J$.}
     \item [(B2)] \textit{Free disposability: $(x,y)\in T$, $(x',y')\in (\Kcp)^{M+S}$, such that  $x'\succeqq x$, $y'\preceqq y $, then $(x',y')\in T$.}
    \item [(B3)] \textit{Convexity: $({x},{y}), (x',y')\in T$, $\lambda\in[0,1]$, then $\lambda (x,y)+(1-\lambda)(x',y') \in T$.}
    \item [(B4)] \textit{Scalability: $(x,y)\in T$, and $\lambda\ge 0$, then $\lambda (x,y)\in  T$.}
\end{itemize}

\begin{theorem}\label{th: characterization fuzzy PPS}
Under axioms (B1), (B2), (B3), and (B4), the CRS interval production possibility set PPS that results from the minimum extrapolation principle are
$$
T_{IDEA} = PPS(X,Y) = \left\{({x},{y})\in (\Kcp)^{M+S} : x \succeqq \sum_{j=1}^{N}\lambda_j {X}_j, 
y \preceqq \sum_{j=1}^{N}\lambda_j {Y}_j, \lambda_j\geq 0, \forall j \right\}
$$
\end{theorem}

\begin{proof}
This is similar to the proof given by Arana-Jimenez et al. \cite{arana2021}.
\end{proof}
In the case of convexity, but not scalability axiom, the PPS is as follows.
\begin{theorem}\label{th: characterization fuzzy PPS convex}
Under axioms (B1), (B2), and (B3), the VRS interval production possibility set PPS that results from the minimum extrapolation principle is
$$
T_{IDEA} = PPS(X,Y) = \left\{({x},{y})\in (\Kcp)^{M+S} : x \succeqq \sum_{j=1}^{N}\lambda_j {X}_j, 
y \preceqq \sum_{j=1}^{N}\lambda_j {Y}_j, \sum_{j=1}^{N}\lambda_j=1, \lambda_j\geq 0, \forall j \right\}
$$
\end{theorem}

\begin{proof}
This is similar to the proof given by Arana-Jimenez et al. \cite{arana2021}.
\end{proof}

\begin{definition}\label{def:efficient}
    A DMU $p$ is said to be efficient if and only if for any $({x},{y})\in T_{IDEA}$, with ${x}\preceqq {X}_p$ and ${y}  \succeqq {Y}_p$ implies $({x},{y})=({X}_p,{Y}_p)$.
\end{definition}{}

% %%%%%%%%%%%%%%%%%%%%%%%%%%%%%%%%%%%%%%%%%%
Given the $T_{IDEA}$, we can recover the following slacks-based measure of the inefficiency interval DEA (IDEA) model with two phases.
%%%%%%%%%%%%
\begin{eqnarray}\label{prob: HIDEA}
\mbox{(IDEA)}\ \ {I}(X_p, Y_p)= 
&\mbox{Max} & {\displaystyle\sum_{i=1}^{M} \frac{\underline{s_i^x}+\overline{s_i^x}}{\underline{x_{ip}}+\overline{x_{ip}}} +
\sum_{r=1}^{S} \frac{\underline{s_r^y}+\overline{s_r^y}}{\underline{y_{rp}}+\overline{y_{rp}}}}
\\%[0.25em]
&\mbox{s.t.} & {\displaystyle\sum_{j=1}^{N}} \lambda_j {x}_{ij} \preceqq {x}_{ip} - s_i^x ,\quad i\in O^{X},
\nonumber\\%[0.25em]
&& {\displaystyle\sum_{j=1}^{N}} \lambda_j {y}_{rj} \succeqq {y}_{rp} + s_r^y ,\quad r\in O^{Y},
\nonumber\\%[0.25em]
&& {\displaystyle\sum_{j=1}^{N}} \lambda_j=1 \nonumber \\%0.25em]
&& \lambda_j \geq 0 ,\quad j=1,\dots, N, \nonumber \\%[0.25em]
&& s_i^x ,s_r^y \in \Kcp, \quad i\in O^{X}, \ \ r\in O^{Y}. \nonumber
\end{eqnarray}
%%%%%%%%%%%%%
\noindent 
%%%%%%
If a $DMU_p$ is efficient, then $I(X_p,Y_p)=0$ (see Arana-Jimenez et al. \cite{arana2021}) but $I(X_p,Y_p)=0$ is insufficient to guarantee a DMU's efficiency. Therefore,  given the optimal solution of \eqref{prob: HIDEA}, $({\bm s^{x^*}},{\bm s^{y^*}},{\bm \lambda^*})$, we proceed with the phase II of the method to exhaust all remaining input and output slacks.
% %%%%%%%
\begin{eqnarray}\label{prob: P2HIDEA}
\mbox{$(PIDEA)_2$}\ \ {H}(X_p, Y_p)= 
&\mbox{Max} & {\displaystyle\sum_{i=1}^{M} \frac{{L_i^x}+{R_i^x}}{\underline{x_{ip}}+\overline{x_{ip}}}}+
\sum_{r=1}^{S} \frac{{L_r^y}+{R_r^y}}{\underline{y_{rp}}+\overline{y_{rp}}}
\\[0.25em]
&\mbox{s.t.} & {\displaystyle\sum_{j=1}^{N}} \lambda_j \underline{x_{ij}} \leq \underline{x_{ip}} - \overline{s_i^{x^*}}  - {R_i^{x}} ,
\quad i\in O^{X},
\nonumber\\
&& {\displaystyle\sum_{j=1}^{N}} \lambda_j \overline{x_{ij}} \leq \overline{x_{ip}} - \underline{s_i^{x^*}} - {L_i^{x}},
\quad i\in O^{X},
\nonumber\\
&& {\displaystyle\sum_{j=1}^{N}} \lambda_j \underline{y_{rj}} \geq \underline{y_{rp}} + \underline{s_r^{y^*}} + {L_r^{y}} ,
\quad  r\in O^{Y},
\nonumber\\
&& {\displaystyle\sum_{j=1}^{N}} \lambda_j \overline{y_{rj}} \geq \overline{y_{rp}} + \overline{s_r^{y^*}} + {R_r^{y}},
\quad  r\in O^{Y},
\nonumber\\
&&{\displaystyle\sum_{j=1}^{N}} \lambda_j=1 \nonumber \\
%&& \lambda_j \geq 0 ,\quad j=1,\dots, N, \nonumber \\[0.25em]
&& \lambda_j , {L_i^x},{R_i^x},{L_r^y},{R_r^y} \ge 0, \quad \forall j, \ \  i\in O^{X}, \ \ r\in O^{Y}. \nonumber
\end{eqnarray}
%%%%%%%%%%%%%%%%%%%%%%%%%%%%%%%%%%%%%%%%%%%%%%%%%
Given a $DMU_p$ with ${I}({X}_p, {Y}_p)=0$, then ${H}({X}_p, {Y}_p)=0$ if and only if $DMU_p$ is efficient. In other words, a $DMU_p$ is efficient if and only if both ${I}({X}_p, {Y}_p)=0$ and ${H}({X}_p, {Y}_p)=0$.
Let $({\bm s^{x^*}},{\bm s^{y^*}},{\bm \lambda^*})$ be the optimal solution of \eqref{prob: HIDEA}, and $\tuple{{\bm L^{x^*}},{\bm R^{x^*}},{ \bm L^{y^*}},{\bm R^{y^*}},{\bm \lambda^{**}}}$ the optimal solution of \eqref{prob: P2HIDEA} for a given $DMU_p$, then we can compute its input and output targets ${X}_{p}^{target}$ and ${Y}_{p}^{target}$ as
\begin{align}
\underline{{x}_{ip}^{target}} &= \underline{{x}_{ip}} - \overline{s_i^{x^*}}-R_i^{x^*}, 
\quad
& \overline{{x}_{ip}^{target}} = \overline{{x}_{ip}} - \underline{s_i^{x^*}}-L_i^{x^*}, 
\qquad i\in O^X, 
\label{eq:xtarget H}
\\[1ex]
\underline{{y}_{rp}^{target}} &=\underline{{y}_{rp}} + \underline{s_r^{y^*}} + L_r^{y^*} ,
\quad
& \overline{{y}_{rp}^{target}} =\overline{{y}_{rp}} + \overline{s_r^{y^*}} + R_r^{y^*} ,
\qquad r\in O^Y. 
\label{eq:ytarget H}
\end{align}
%%%%%%%%%%
\noindent 
As discussed in \cite{arana2021}, considering new slacks $L_i^x,L_r^y$, and $R_i^x,L_r^y$ in $(PIDEA)_2$, or phase 2, was necessary to guarantee the efficiency of the corresponding DMU, as well a to provide their targets. This is because the feasible slacks improvements of the inequality constraints in (\ref{prob: HIDEA}) are not exhausted. 
Given two positive intervals $A,B \in \Kcp$, with $A\preceqq B$, the positive interval upper slack $su\in \Kcp$ such that  $A = B-su$, or the lower slack $sl \in \Kcp$ such that  $A+sl = B$ do not necessarily exist. 
For instance, if we consider $A=[2,3]$ and $B=[2,5]$, then 
$A=[2,3] = B-su = [2-\overline{su},5-\underline{su}]$.  
It is not possible to find a positive interval upper slack $su$ holding the equality. 
But from $A + sl =[2 + \underline{sl}, 3 + \overline{sl}] = [2,5] = B $, we can find $sl = [0,2]$.
Similarly, for $A=[2,6]$ and $B=[5,7]$, the positive interval lower slack $sl \in \Kcp$ such that  $A+sl = B$ does not necessarily exist. But we find $su = [1,3]$ such that $A = B - su$.
This is, given $A$ and $B$ two non negative closed intervals, with $A\preceqq B$, there always exists $su$ or $sl$ in $\Kcp$ such that $A = B-su$ or $A+sl = B$, as we proof following.
%%%%%%%%%
\begin{proposition}\label{pro:slacks}
If $A, B\in\Kcp$, with $A\preceqq B$, then there exist $sl, su\in\Kcp$ such that $A = B-su$ or $A+sl = B$. 
\end{proposition}
%%%%%%%%%%
\begin{proof}
On one hand, following Stefanini and Arana \cite{stefanini-arana19}, if $A\preceqq B$, then $B\ominus_{gH}A\succeqq 0$. On the other hand, there always exists $C\in\Kcp$ such that
$$B\ominus_{gH}A=C \Longleftrightarrow
\left\{
\begin{array}{ll}
&(a)\ B= A+C,\\
or&(b)\ A = B + (-1)C.
\end{array}
\right.$$
Since $C=B\ominus_{gH}A\succeqq 0$, then $C\in\Kcp$. In case (a), we have that $B= A+C$, we define $sl=C$; consequently, the proposition holds. And in case (b), that is, $A = B + (-1)C$, we define $su=C$, and the proof is complete.
\end{proof}
\begin{proposition}\label{pro:slacks 2}
If $A, B\in\Kcp$, and $sl, su\in\Kcp$ such that $A +sl= B-su$, then $A\preceqq B$. 
\end{proposition}
%%%%%%%%%%
\begin{proof}
Since $A\preceqq A +sl$, and $B-su\preceqq B$, then it follows that $A\preceqq A +sl=B-su\preceqq B$.
\end{proof}
Based on the previous propositions, we get the following useful result for formulating the forthcoming models.
\begin{corollary}\label{cor:slacks}
Given $A, B\in\Kcp$, then $A\preceqq B$ if and only if then there exist $sl, su\in\Kcp$ such that $A+sl=B-su$, with $sl=0$ or $su=0$.
\end{corollary}
\begin{proof}
It was shown in the proof of Proposition \ref{pro:slacks} that, if $A\preceqq B$, there exists $sl\in\Kcp$ or $su\in\Kcp$, depending on the case (a) or (b), and we then have $su= 0 $ or $sl= 0$ respectively. In the reverse direction, if $A+sl=B-su$, with $sl, su\in\Kcp$, then, from Proposition \ref{pro:slacks 2}, it follows that $A\preceqq B$.
\end{proof}

\begin{remark}\label{r: 1}
Given the definition of the {$gH$-difference}, Eq.\eqref{eq:gH diff bis}, Corollary \ref{cor:slacks} actually implies that $sl = B\ominus_{gH}A $ and $su = 0$, when $\underline{b}-\underline{a} \leq \overline{b}-\overline{a}$ (i.e., 
$A + sl = [\underline{a} , \overline{a}] + 
[ \underline{b}-\underline{a} \ , \ \overline{b}-\overline{a} ]
= [\underline{b}, \overline{b}]  = 
B - su $ ), and  
 $sl =0$ and $su =  B\ominus_{gH}A $,  otherwise (i.e., 
$ A + sl = [\underline{a}, \overline{a}] =  
[\underline{b} , \overline{b}] - [\overline{b}-\overline{a} \ , \ \underline{b}-\underline{a}]  = B - su $). 
\end{remark}

\begin{remark}\label{r: 2}
In the previous result, if $A$ and $B$ are crisp, with $sl=0$ or $su=0$, it implies that both $sl$ and $su$ are crisp.
\end{remark}

In the next section, we present our proposed model based on applying these previous results to the (IDEA) model (\ref{prob: HIDEA}), aiming to unify the two phases (IDEA and PIDEA) and exhausting the constraints to equalities at once.

%%%%%%%%%%%%%%%%%%%%%%%%%%%%%%%%%%%%
\section{Enhanced interval slack-based model}
%%%%%%%%%%%%%%%%%%%%%%%%%%%%%%%%%%%%

Let us recover the (IDEA) model \eqref{prob: HIDEA}, and take into account the previous discussions on its constraints, as well as Proposition \ref{pro:slacks}. Thus, we propose the following Enhanced Inefficiency Non-Linear program (EINL) for our current DEA framework with interval data .
%%%%%%%
\begin{eqnarray}%\label{prob: EINL}
\mbox{(EINL)}\ \ {EI}(X_p, Y_p)= 
&\mbox{Max} & {\displaystyle\sum_{i=1}^{M} \frac{\underline{sl_i^x}+\underline{su_i^x}+\overline{sl_i^x}+\overline{su_i^x}}{\underline{x_{ip}}+\overline{x_{ip}}} +
\sum_{r=1}^{S} \frac{\underline{sl_r^y}+\underline{su_r^y}+\overline{sl_r^y}+\overline{su_r^y}}{\underline{y_{rp}}+\overline{y_{rp}}}}
\hspace{1.5cm}
\label{prob: EINL}\\[0.25em]
&\mbox{s.t.} & {\displaystyle\sum_{j=1}^{N}} \lambda_j {x}_{ij}+ sl_i^x = {x}_{ip} - su_i^x ,\quad i=1,\ldots,M,
\label{prob: EINL1}\\[0.25em]
&& {\displaystyle\sum_{j=1}^{N}} \lambda_j {y}_{rj}-su_r^y =
{y}_{rp} + sl_r^y ,\quad r=1,\ldots,S,
\label{prob: EINL2}\\[0.25em]
&& {\displaystyle\sum_{j=1}^{N}} \lambda_j=1 \label{prob: EINL3} \\[0.25em]
&& sl_i^x\cdot su_i^x=0 ,\; sl_r^y\cdot su_r^y=0,\quad i=1,\ldots,M, \ r=1,\ldots,S,
\label{prob: EINL4}\\[0.25em]
&& \lambda_j \geq 0 ,\quad j=1,\dots, N, \label{prob: EINL5} \\[0.25em]
&& sl_i^x, su_i^x ,sl_r^y, su_r^y \in \Kcp, \quad i=1,\ldots,M, \ r=1,\ldots,S. \label{prob: EINL6}
\end{eqnarray}
%%%%%%%%%%%%%
\noindent 
Comparing (IDEA) and (EINL), we tighten the input and output constraints to equality, by considering low and upper slack variables on both sides  of the constraints such that only one of them can be non-zero (as imposed in Eq. \eqref{prob: EINL4}). Moreover, as stated in the following result, the inefficiency measure of any DMU under the (EINL) model is always higher than or equal to the inefficiency measure under (IDEA).
%%%%%%%%
\begin{proposition}\label{pro: comparison}
Given any $DMU_p$, it is verified that ${EI}(X_p, Y_p)\ge {I}(X_p, Y_p)$.
\end{proposition}
\begin{proof}
   For a given $DMU_p$, let $({\bm s^{x^*}},{\bm s^{y^*}},{\bm \lambda^*})$ be the optimal solution of (IDEA) \eqref{prob: HIDEA}, and $\tuple{{\bm L^{x^*}},{\bm R^{x^*}},{ \bm L^{y^*}},{\bm R^{y^*}},{\bm \lambda^{**}}}$ the optimal solution of \eqref{prob: P2HIDEA}. 
    If we define  
    ${\bm su^{x}} = \left[ \underline{s_i^{x^*}}+L_i^{x^*} ,  \overline{s_i^{x^*}} + R_i^{x^*} \right]$, 
    ${\bm sl^{x}}=0$, 
    ${\bm sl^{y}} = \left[ \underline{s_r^{y^*}} + L_r^{y^*} , \overline{s_r^{y^*}} + R_r^{y^*} , \right]$, and 
    ${\bm su^{y}}=0$, 
    from \eqref{eq:xtarget H} and \eqref{eq:ytarget H}, it is straightforward that $({\bm sl^{x}},{\bm su^{x}},{\bm sl^{y}},{\bm su^{y}},{\bm \lambda^{**}})$ is a feasible solution for (EINL), and from \eqref{prob: HIDEA} and \eqref{prob: EINL} that ${EI}(X_p, Y_p)\ge {I}(X_p, Y_p)$.
\end{proof}
%%

%%%%%%%%%%%%%%%%%%%%%%%%%%%%%%%%%%%%%%%%%%%%%%%%%
Given the (EINL) model \eqref{prob: EINL}, it is expected that efficient DMUs have a null inefficiency measure. We prove that only efficient DMUs satisfy the latter and, hence, that ${EI}({X}_p, {Y}_p)=0$ is a characterization of efficient DMUs.

\begin{theorem}\label{th: efficiency}
$DMU_p$ is efficient if and only if $EI(X_p,Y_p)=0$.
\end{theorem}
\begin{proof}
By contradiction. Let us assume that ${EI}({X}_p, {Y}_p)=0$. 
Then, if $DMU_p$ is not efficient there exist $({x}^*,{y}^*)\in T_{IDEA}$ such that, ${x}^* \precneqq {X}_p$ and  ${Y}_p \precneqq {y}^*$. %
As $({x}^*,{y}^*)\in T_{IDEA}$, then it also holds that 
${x}^* \succeqq {\sum_{j=1}^{N}} \lambda^*_j {x}_{ij} $ and ${y}^* \preceqq {\sum_{j=1}^{N}} \lambda^*_j {y}_{rj}$ for some $\bm \lambda^* \in \R^N_{+}$, with $\sum_{j=1}^N \lambda^*_j = 1$. 
Then, we have that ${\sum_{j=1}^{N}} \lambda^*_j {x}_{ij} \precneqq X_p$, and 
$  Y_p \precneqq {\sum_{j=1}^{N}} \lambda^*_j {y}_{rj}$.
Applying Corollary \ref{cor:slacks}, we find some 
$sl^{x^{*}},su^{x^{*}}\succeqq 0 $, with 
$sl^{x^{*}}\cdot su^{x^{*}}= 0 $,  such that
$\sum_{j=1}^{N} \lambda^*_j {x}_{ij} + sl_{i}^{x^{*}} = {x}_{ip} - su_{i}^{x^{*}}$, 
and $sl^{x^{*}} \ne 0$, or $su^{x^{*}}\ne 0$. This is due to the $\precneqq$ constraint derived from the definition of efficiency. 
Similarly, in the case of the outputs, there are some  $sl^{y^{*}},su^{y^{*}}\succeqq 0$, with $sl^{y^{*}}\cdot su^{y^{*}}= 0$, such that 
$\sum_{j=1}^{N} \lambda^*_j {y}_{rj} - su_r^{y^{*}} = {Y}_{rp} + sl_r^{y^{*}}$, and $sl^{y^{*}} \ne 0$, or $su^{y^{*}}\ne 0$.
By construction, it is straightforward that 
$\tuple{\bm \lambda^*, sl^{x^{*}}}, su^{x^{*}},
sl^{y^{*}},su^{y^{*}}$ 
is a feasible solution of model $(EINL)$, but its objective function \eqref{prob: EINL} is strictly greater than zero, which is a contradiction with $EI(X_p,Y_p) = 0$.   

Now let us suppose $DMU_p$ is efficient, but ${EI}({X}_p, {Y}_p)>0$. Then there are some 
${sl}_{i_o}^{x}\succeqq 0$ and ${sl}_{i_o}^{x}\ne 0$, or  ${su}_{i_o}^{x} \succeqq 0$ and  ${su}_{i_o}^{x} \ne 0$, for some $i_o \in \{1,\ldots,M\}$. Or there are some  ${sl}_{r_o}^{y} \succeqq 0, \ {sl}_{r_o}^{y} \ne 0$, or ${su}_{r_o}^{y} \succeqq0, \ {su}_{r_o}^{y} \ne 0$, 
for some $r_o \in \{1,\ldots,S\}$. 
We can then compute $(x^*,y^*)$ as 
${x}_i^* = x_{ip}$ for all $i\in \{1,\ldots,M\}$, $i\ne i_o$, and 
${x}_{i_o}^* = x_{i_op} - {su}_{i_o}^{x}$ if ${su}_{i_o}^{x} \ne 0$, 
or ${x}_{i_o}^* ={\sum_{j=1}^{N}} \lambda_j {x}_{i_oj}$ if ${sl}_{i_o}^{x} \ne 0$. 
Analogously, ${y}_r^* = y_{rp}$ for all $r\in \{1,\ldots,S\}$, with $r\ne r_o$, and 
${y}_{r_o}^* = y_{r_op} + {sl}_{r_o}^{y}$ if ${sl}_{r_o}^{y}\ne 0$, 
or ${y}_{r_o}^* ={\sum_{j=1}^{N}} \lambda_j {y}_{r_oj}$ if ${su}_{r_o}^{y} \ne 0$. By definition, $(x^*,y^*)$ holds constraints \eqref{prob: EINL1} and \eqref{prob: EINL2}, then it belongs to the PPS, $({x}^*,{y}^*)\in T$, and 
${x}^*\preceqq {X}_p$ and ${y}^*  \succeqq {Y}_p$, with $({x}^*,{y}^*) \ne ({X}_p,{Y}_p)$. This implies a contradiction with the fact that $DMU_p$ is efficient.

\end{proof}

The previous mathematical program $(EINL)$ is nonlinear, but we can provide an equivalent linear program with $0-1$ variables leading to the following Enhanced Inefficiency Mix-Integer Linear program  (EIMIL), 

\begin{eqnarray}%\
\mbox{(EIMIL)}\ \ {EI}(X_p, Y_p)= 
&\mbox{Max} & {\displaystyle\sum_{i=1}^{M} \frac{\underline{sl_i^x}+\underline{su_i^x}+\overline{sl_i^x}+\overline{su_i^x}}{\underline{x_{ip}}+\overline{x_{ip}}} +
\sum_{r=1}^{S} \frac{\underline{sl_r^y}+\underline{su_r^y}+\overline{sl_r^y}+\overline{su_r^y}}{\underline{y_{rp}}+\overline{y_{rp}}}}
\hspace{1cm}
\label{prob: EIMIL}\\[0.25em]
&\mbox{s.t.} & \eqref{prob: EINL1}-\eqref{prob: EINL3},\eqref{prob: EINL5}-\eqref{prob: EINL6} 
\nonumber\\[0.25em]
&& sl_i^x \preceqq  L^x_{i} z^x_{i},\quad i=1,\ldots,M,
\label{prob: EIMIL1}\\[0.25em]
&& su_i^x\preceqq  R^x_{i} (1-z^x_{i}),\quad i=1,\ldots,M,
\label{prob: EIMIL2}\\[0.25em]
&& sl_r^y \preceqq  L^y_{r} z^y_{r},\quad r=1,\ldots,S,
\label{prob: EIMIL3}\\[0.25em]
&& su_r^y\preceqq  R^y_{r} (1-z^y_{r}),\quad r=1,\ldots,S,
\label{prob: EIMIL4}\\[0.25em]
&& z^x_i, z^y_r \in  \{ 0,1\},\quad i=1,\ldots,M, \ r=1,\ldots,S.
\label{prob: EIMIL5}
\end{eqnarray}
%%%%%%%%%%%%%
\noindent 
where equations \eqref{prob: EIMIL1} to \eqref{prob: EIMIL5} are equivalent to the non-linear constrains \eqref{prob: EINL4}, and $ L^x_{i},R^x_{i},L^y_{r},R^y_{r}$ are real positive constants, and large enough. Recall that a number can be identified with an interval whose extremes are equal and coincide with such a number. Then, it can be compared with intervals via interval inequalities. 
We can consider 
$L^x_{i}=R^x_{i}=\overline{x_{ip}}$, and $L^y_{r}=R^y_{r}=\max\{ \overline{y_{rp}}: r=1,\ldots,S \}$. %Then, 
Another possibility is to set $L^x_{i}=R^x_{i}=L^y_{r}=R^y_{r}=\max\{ \overline{x_{ip}},\overline{y_{rp}}: i=1,\ldots,M, r=1,\ldots,S \}$.
Thus, (EIMIL) is an equivalent mix-integer linear program formulation to (EINL), which computes the same Enhanced Inefficiency measure ${EI}({X}_p, {Y}_p)$.

Let $({\bm sl^{x^*}}, {\bm su^{x^*}},{\bm sl^{y^*}}, {\bm sl^{y^*}},{\bm \lambda^*}, {\bm z^{x^*}},{\bm z^{y^*}})$ be an optimal solution for $(EIMIL)$ model (\ref{prob: EIMIL}) for a given $DMU_p$, then we can compute its input and output targets 
${X}_{p}^{target} = \bm \lambda^* \cdot X = \sum_{j = 1}^N \lambda^*_j X_j$  
and 
${Y}_{p}^{target} = \bm \lambda^* \cdot Y = \sum_{j = 1}^N \lambda^*_j Y_j$ as
%%%%%
\begin{align}
\underline{{x}_{ip}}^{target} &= \underline{{x}_{ip}}  -\overline{ su_i^{x^*}}  - \underline{sl_i^{x^*}} , 
&\overline{{x}_{ip}}^{target} &= \overline{{x}_{ip}} - \underline{su_i^{x^*}}-\overline{sl_i^{x^*}}  , \ \  i=1,\ldots,M,
\label{eq:xtarget}
\\[1.ex]
\underline{{y}_{rp}}^{target} & = \underline{{y}_{rp}} + \underline{sl_r^{y^*}} + \overline{su_r^{y^*}} ,
& \overline{{y}_{rp}}^{target} &= \overline{{y}_{rp}} + \overline{sl_r^{y^*}} + \underline{su_r^{y^*}}, \ \   r=1,\ldots,S.
\label{eq:ytarget}
%\\ \nonumber
\end{align}
The above Eq. \ref{eq:xtarget} and \ref{eq:ytarget} are just the parametrizations derived from 
$X_p^{target} + sl^{x^* } = X_p - su^{x^*}$ \eqref{prob: EINL1}, 
and $Y_p^{target} - su^{y^*} = Y_p + sl^{y^*}$ \eqref{prob: EINL2}, respectively.

\begin{theorem}\label{th: efficiency targets}
$({X}_{p}^{target},{Y}_{p}^{target})$ is efficient.
\end{theorem}

\begin{proof}
Let $({\bm sl^{x^*}}, {\bm su^{x^*}},{\bm sl^{y^*}}, {\bm sl^{y^*}},{\bm \lambda^*}, {\bm z^{x^*}},{\bm z^{y^*}})$ be an optimal solution for $(EIMIL)$ model (\ref{prob: EIMIL}) for a given $DMU_p$. 
From Eq. \eqref{prob: EINL1}-\eqref{prob: EINL2} and \eqref{eq:xtarget}-\eqref{eq:ytarget}, it follows that $({X}_{p}^{target},{Y}_{p}^{target}) \in T_{IDEA}$, and ${X}_{p}^{target} \preceqq X_p$, ${Y}_{p}^{target} \succeqq Y_p$ (Proposition \ref{pro:slacks 2}).

If $({X}_{p}^{target},{Y}_{p}^{target})$ were not efficient there would be some $(x',y')\in T_{IDEA}$ such that ${x}'\preceqq {X}_{p}^{target}$, $ {Y}_{p}^{target} \preceqq {y}'$, with
${x}'\precneqq {X}_{p}^{target}$ or $ {Y}_{p}^{target} \precneqq {y}'$, and  
$x' \succeqq \sum_{j=1}^{N}\lambda'_j X_j$,
$ y' \preceqq \sum_{j=1}^{N}\lambda'_j Y_j$, for some $\bm \lambda' \in \R^N_{+}$, with $\sum_{j=1}^N \lambda_j = 1$.
In summary, we have that
$$
\sum_{j=1}^{N}\lambda'_j X_j  \preceqq  
    x'  \precneqq X^{target}_p 
     \preceqq   X_p, \quad {or} \quad
    Y_p   \preceqq  Y^{target}_p 
    \precneqq y'   \preceqq 
    \sum_{j=1}^{N}\lambda'_j Y_j
$$
In the case of the inputs, applying Corollary \ref{cor:slacks}, we can find some slacks $0\precneqq  sl^{x'},  su^{x'}\in \Kcp $, such that 
$sl^{x'}\cdot su^{x'} = 0$, and $ \sum_{j=1}^{N}\lambda'_j X_j - su^{y'} = X_p +  sl^{y'}$.
Analogously, for the output case, we can find some slacks $0\precneqq sl^{y'}, sl^{y'}\in \Kcp $, such that 
$sl^{y'}\cdot su^{y'} = 0$, and $ \sum_{j=1}^{N}\lambda'_j Y_j - su^{y'} = Y_p +  sl^{y'}$.
Moreover, given Remark \ref{r: 1} and that all intervals are $\Kcp$, it follows 
$$
\left.
\begin{array}{l}
   sl_i^{x^*} = x_{ip} \ominus_{gH}  \sum_{j=1}^{N}\lambda^*_j x_{ij}, \ \   su_i^{x^*} = 0, %\quad \text{or}  
   \\[1ex]
   sl_i^{x^*} = 0, \ \  su_i^{x^*} = x_{ip} \ominus_{gH}  \sum_{j=1}^{N}\lambda^*_j x_{ij}
\end{array} \right\}
, \text{ and} \quad 
\left.
\begin{array}{l}
   sl_i^{x'} = x_{ip}\ominus_{gH}  \sum_{j=1}^{N}\lambda'_j x_{ij}, \ \ su_i^{x'} = 0,
   \\[1ex]
   sl_i^{x'} = 0, \ \ su_i^{x'} = x_{ip}\ominus_{gH}  \sum_{j=1}^{N}\lambda'_j x_{ij}
\end{array} \right\}
$$
$$\Rightarrow \quad 
sl_i^{x^*}, \text{ or } su_i^{x^*} = x_{ip} \ominus_{gH}  \sum_{j=1}^{N}\lambda^*_j x_{ij} 
\ \ \precneqq \ \ 
x_{ip}\ominus_{gH}  \sum_{j=1}^{N}\lambda'_j x_{ij} = 
sl_i^{x'}, \text{ or } su_i^{x'}.
$$
Analogously, for the output case, it holds  
$$
\left.
\begin{array}{l}
   sl_r^{y^*} = \sum_{j=1}^{N}\lambda^*_j y_{rj} \ominus_{gH}  y_{rp} , \ \   su_r^{y^*} = 0, %\quad \text{or}  
   \\[1ex]
   sl_r^{y^*} = 0 , \ \   su_r^{y^*} = \sum_{j=1}^{N}\lambda^*_j y_{rj} \ominus_{gH}  y_{rp},
\end{array} \right\}
, \text{ and} \quad 
\left.
\begin{array}{l}
    sl_r^{y'} = \sum_{j=1}^{N}\lambda'_j y_{rj} \ominus_{gH}  y_{rp} , \ \   su_r^{y'} = 0, %\quad \text{or}  
   \\[1ex]
   sl_r^{y'} = 0 , \ \   su_r^{y'} = \sum_{j=1}^{N}\lambda'_j y_{rj} \ominus_{gH}  y_{rp},
\end{array} \right\}
$$
$$\Rightarrow \quad 
sl_r^{y^*}, \text{ or } su_r^{y^*} =\sum_{j=1}^{N}\lambda^*_j y_{rj} \ominus_{gH}  y_{rp} 
\ \ \precneqq \ \ 
\sum_{j=1}^{N}\lambda'_j y_{rj} \ominus_{gH}  y_{rp}  = 
sl_r^{y'}, \text{ or } su_r^{y'}.
$$

Therefore we can find a feasible solution for model $(EIMIL)$, $\tuple{ \bm su^{x'}, \bm sl^{x'}, \bm  sl^{y'}, \bm su^{y'}, {\bm \lambda'}, {\bm z^{x'}},{\bm z^{y'}}}$, with a larger objective function value \eqref{prob: EIMIL} than the optimum, $({\bm sl^{x^*}}, {\bm su^{x^*}},{\bm sl^{y^*}}, {\bm sl^{y^*}},{\bm \lambda^*}, {\bm z^{x^*}},{\bm z^{y^*}})$, which is a contradiction. 

\end{proof}
\noindent
To solve the (EIMIL) model \eqref{prob: EIMIL}, we can reformulate it as
%%%%%%%%
\begin{eqnarray}%\label{prob: ELPIDEA}
\mbox{(PEIMIL)}\ \ {EI}(X_p, Y_p)= 
&\mbox{Max} & {\displaystyle\sum_{i=1}^{M} \frac{\underline{sl_i^x}+\underline{su_i^x}+\overline{sl_i^x}+\overline{su_i^x}}{\underline{x_{ip}}+\overline{x_{ip}}} +
\sum_{r=1}^{S} \frac{\underline{sl_r^y}+\underline{su_r^y}+\overline{sl_r^y}+\overline{su_r^y}}{\underline{y_{rp}}+\overline{y_{rp}}}}
\hspace{1.cm} 
\label{prob: ELPIDEA} 
\\[0.25em]
&\mbox{s.t.} & {\sum_{j=1}^{N}} \lambda_j \underline{{x}_{ij}}+ \underline{sl_i^x} = \underline{{x}_{ip}}  -\overline{ su_i^x}  ,\quad i=1,\ldots,M,
\label{prob: ELPIDEA1}
\\%[0.25em]
&& {\sum_{j=1}^{N}} \lambda_j \overline{{x}_{ij}} +\overline{ sl_i^x}  =  \overline{ {x}_{ip}}  -\underline{ su_i^x}  ,\quad i=1,\ldots,M,
\label{prob: ELPIDEA2}
\\%[0.25em]
&& {\sum_{j=1}^{N}} \lambda_j \underline{{y}_{rj}}- \overline{ su_r^y} =  \underline{{y}_{rp}} +\underline{ sl_r^y} ,\quad r=1,\ldots,S,
\label{prob: ELPIDEA3}
\\%[0.25em]
&& {\sum_{j=1}^{N}} \lambda_j \overline{{y}_{rj}}-\underline{su_r^y} =  \overline{{y}_{rp}} +\overline{ sl_r^y} ,\quad r=1,\ldots,S,
\label{prob: ELPIDEA4}
\\%[0.25em]
&& {\sum_{j=1}^{N}} \lambda_j=1 , 
\label{prob: ELPIDEA5} 
\\%[0.25em]
&& \lambda_j \geq 0 ,\quad j=1,\dots, N, 
\label{prob: ELPIDEA6} 
\\%[0.25em]
&& \underline{sl_i^x} \leq  \overline{sl_i^x}\leq  L^x_{i} z^x_{i}, 
\quad i=1,\ldots,M,
\label{prob: ELPIDEA7}
\\%[0.25em]
&&\underline{ su_i^x}   \leq \overline{su_i^x}\leq  R^x_{i} (1-z^x_{i}),  
\quad i=1,\ldots,M,
\label{prob: ELPIDEA8}
\\%[0.25em]
&& \underline{sl_r^y} \leq  \overline{ sl_r^y} \leq  L^y_{r} z^y_{r},
\quad r=1,\ldots,S,  
\label{prob: ELPIDEA9}
\\%[0.25em]
&& \underline{ su_r^y} \leq  \overline{  su_r^y}\leq  R^y_{r} (1-z^y_{r}),
\quad r=1,\ldots,S,
\label{prob: ELPIDEA10}
\\%[0.25em]
&&\underline{ sl_i^x} , \overline{ sl_i^x}, \underline{ su_i^x} , \overline{ su_i^x} \geq  0, \quad i=1,\ldots,M,
\label{prob: ELPIDEA11}
\\%[0.25em]
&&  \underline{ sl_r^y}, \overline{ sl_r^y},  \underline{ su_r^y}, \overline{ su_r^y}\geq  0, \quad r=1,\ldots,S,
\label{prob: ELPIDEA12}
\\%[0.25em]
&& z^x_i, z^y_r \in  \{ 0,1\},\quad i=1,\ldots,M, \ r=1,\ldots,S.
\label{prob: ELPIDEA13}
\end{eqnarray}
%%%%%%%%%%%%%

%%%%%%%%%%%%%%%%%%%%%%%%%%%%%%%%%%%%%%%%%%%%%%%%%
\section{Super SBI interval model}

The above (EIMIL) model results indicate if DMUs are efficient, but all the efficient DMUs obtain the same zero inefficient measure. Hence, to discriminate between efficient units, a super-efficiency approach is considered. Super-efficiency as originally proposed by Andersen and Petersen \cite{Andersen01} involved radial DEA models. The super-efficiency method applied to the Slack-based measure of inefficiency (SBI) metrics is appropriately called the super SBI model (see, e.g., Moreno and Lozano \cite{Moreno01}). The idea behind the super-efficiency concept is to exclude the observation being benchmarked from the set of observations that define the technology. 
We present the proposed non-oriented super SBI DEA model, which should be solved only for DMUs labeled efficient by the (EIMIL) model.

\begin{eqnarray} %\label{prob: SEIMIL}
\mbox{(S-EIMIL)}\ \ {SEI}(X_p, Y_p)= 
&\mbox{Min} & {\displaystyle\sum_{i=1}^{M} \frac{\underline{sl_i^x}+\underline{su_i^x}+\overline{sl_i^x}+\overline{su_i^x}}{\underline{x_{ip}}+\overline{x_{ip}}} +
\sum_{r=1}^{S} \frac{\underline{sl_r^y}+\underline{su_r^y}+\overline{sl_r^y}+\overline{su_r^y}}{\underline{y_{rp}}+\overline{y_{rp}}}}
\hspace{1.cm}
\label{prob: SEIMIL}
\\[0.25em]
&\mbox{s.t.} & {\displaystyle\sum_{j=1, j\ne p}^{N}} \lambda_j {x}_{ij} {-} sl_i^x \preceqq {x}_{ip} {+} su_i^x ,\quad i=1,\ldots,M,
\label{prob: SEIMIL1}\\[0.25em]
&& {\displaystyle\sum_{j=1, j\ne p}^{N}} \lambda_j {y}_{rj} {+} su_r^y \succeqq {y}_{rp} {-} sl_r^y ,\quad r=1,\ldots,S,
\label{prob: SEIMIL2}\\[0.25em]
% && \lambda_j \geq 0 ,\quad j=1,\dots, N, \nonumber \\[0.25em]
&& {\displaystyle\sum_{j=1, j\ne p}^{N}} \lambda_j=1 
\label{prob: SEIMIL3}\\[0.25em]
&& \lambda_j \geq 0, \ \ j\ne p,
\label{prob: SEIMIL4}\\[0.25em]
&& \eqref{prob: EINL6} , \eqref{prob: EIMIL1}-\eqref{prob: EIMIL5}
\nonumber
\end{eqnarray}
%%%%%%%%%%%%%
where $ L^x_{i}, R^x_{i}, L^y_{r}, R^y_{r}$ are real positive constant numbers and great enough, as commented before.
The corresponding equivalent parametrized model is: 

\begin{eqnarray}\label{prob: SPEIMIL}
\mbox{(S-PEIMIL)}\ \ {SEI}(X_p, Y_p)= 
&\mbox{Min} & {\displaystyle\sum_{i=1}^{M} \frac{\underline{sl_i^x}+\underline{su_i^x}+\overline{sl_i^x}+\overline{su_i^x}}{\underline{x_{ip}}+\overline{x_{ip}}} +
\sum_{r=1}^{S} \frac{\underline{sl_r^y}+\underline{su_r^y}+\overline{sl_r^y}+\overline{su_r^y}}{\underline{y_{rp}}+\overline{y_{rp}}}}
\hspace{1.cm}
\label{prob: SPEIMIL}\\[0.25em]
&\mbox{s.t.} & {\displaystyle\sum_{j=1,j\ne p}^{N}} \lambda_j \underline{{x}_{ij}} - \overline{sl_i^x} \leq \underline{{x}_{ip}}  + \underline{ su_i^x}  ,\quad i=1,\ldots,M,
\label{prob: SPEIMIL1}\\[0.25em]
&& {\displaystyle\sum_{j=1,j\ne p}^{N}} \lambda_j \overline{{x}_{ij}}  - \underline{ sl_i^x}  \leq \overline{ {x}_{ip}}  + \overline{ su_i^x}  ,\quad i=1,\ldots,M,
\label{prob: SPEIMIL2}\\[0.25em]
&& {\displaystyle\sum_{j=1,j\ne p}^{N}} \lambda_j \underline{{y}_{rj}} + \underline{ su_r^y} \geq \underline{{y}_{rp}} - \overline{ sl_r^y} ,\quad r=1,\ldots,S,
\label{prob: SPEIMIL3}\\[0.25em]
&& {\displaystyle\sum_{j=1,j\ne p}^{N}} \lambda_j \overline{{y}_{rj}} + \overline{su_r^y} \geq \overline{{y}_{rp}} - \underline{ sl_r^y} ,\quad r=1,\ldots,S,
\label{prob: SPEIMIL4}\\[0.25em]
% && {\displaystyle\sum_{j=1,j\ne p}^{N}} \lambda_j=1 , \label{prob: SPEIMIL5} \\[0.25em]
&& \eqref{prob: ELPIDEA7}-\eqref{prob: ELPIDEA13},  \eqref{prob: SEIMIL3}-\eqref{prob: SEIMIL4}.
\nonumber
\end{eqnarray}
%%%%%%%%%%%%%

Let us note that in the crisp case, that is, when all data are crisp, then, from Remark \ref{r: 2}, all slacks in (SuperEIMIL) are crisp, and all interval equalities and inequalities become ordinary equalities and inequalities. Hence, in that case, defining $s_i^x=sl_i^x+ su_i^x$, $s_r^y=sl_r^y+ su_r^y$, the problem reduces to the conventional super SBI DEA model, i.e.

\begin{eqnarray}\label{prob: superDEA}
\mbox{(SuperDEA)}\ \ {SuperI}(X_p, Y_p)= 
&\mbox{Min} & {\displaystyle\sum_{i=1}^{M} \frac{s_i^x}{x_{ip}} +
\sum_{r=1}^{S} \frac{s_r^y}{y_{rp}}}
\\[0.25em]
&\mbox{s.t.} & {\displaystyle\sum_{j=1, j\ne p}^{N}} \lambda_j {x}_{ij} \leq {x}_{ip} + s_i^x ,\quad i=1,\dots, M,
\nonumber\\[0.25em]
&& {\displaystyle\sum_{j=1, j\ne p}^{N}} \lambda_j {y}_{rj} \geq {y}_{rp} - s_r^y ,\quad r=1,\dots, S,
\nonumber\\[0.25em]
&& {\displaystyle\sum_{j=1, j\ne p}^{N}} \lambda_j=1 \nonumber \\[0.25em]
&& \lambda_j \geq 0 ,\quad j=1,\dots, N, \nonumber \\[0.25em]
&& s_i^x ,s_r^y \geq 0 , \quad i=1,\dots, M, \ \ r=1,\dots, S. \nonumber
\end{eqnarray}

%%%%%%%%%%%%%%%%%%%%%%%%%
\section{Application to Tourism}\label{sec:Example}

This section aims to assess, using data envelopment analysis (DEA), the sustainability efficiency of tourism in the most important Mediterranean regions during 2019. According to the World Tourism Organization, sustainable development is “tourism that takes full account of its current and future economic, social, and environmental impacts, addressing the needs of visitors, the industry, the environment, and host communities.” Sustainability is usually represented in three fundamental pillars or dimensions: economic, social, and environmental \cite{Lozano7}. Sustainability is a recent concept that is very important nowadays for the following reasons:

\begin{itemize}
\item A key to preserving the planet
\item It helps to reduce pollution and conserve resources
\item Creating jobs and stimulating the economy 
\item Improves public health
\item Protects biodiversity 
\item A development that is achievable with political will and public support
\end{itemize}

On the other hand, tourism is considered one of the leading international commerce sectors and one of the primary sources of income for many developing countries. During the last decades, tourism has represented a vital world business and has experienced continued growth. For example, according to World Tourism Organization, international tourism arrivals grew 4,3\% in 2014, reaching 1.133 million tourists, and in January-March 2019 compared to the same period last year, below the 6\% average growth of the past two years. However, it is essential to note that 2020 was a challenging year for most sectors cause of the Covid-19 pandemic, and the tourism industry was affected notably. The number of tourist trips undertaken each year before the advent of Covid-19 exceeded the world's population \cite{Rasoolimanesh}. Although according to the latest World Tourism Barometer from UNWTO, this is in the end because international tourism is on track to reach 65\% of pre-pandemic levels by the end of 2022, and the sector continues to recover from the pandemic.

Recent studies on sustainable tourism using Data Envelopment Analysis focus on the environmental effects and competitiveness of tourism. For example, Huang et al. \cite{Huang01} use SBM-DEA and Tobit's regression to measure the efficiency of environmental training for diving tourists considering inputs such as education, Diver's qualifications, or length of diving time and output as improper environmental behaviors. Also, regarding eco-efficiency, Li et al. \cite{Li01} use two DEA models (CCR and Panel Tobit) to assess the Chinese forest parks in 30 provinces of China, considering inputs as forest park employees, ecological tourism footprint, water consumption or annual forest park tourism data and outputs as total tourism revenue, SO2 emissions or solid particulate emissions. In the matter of evaluating the impact of high-speed rail on the development efficiency of low-carbon in China, Li et al. \cite{Li02} considering an input (namely, high-speed rail) and an output (namely, low-carbon tourism) through the stochastic production frontier method (SFA) in combination with BCC-DEA models. Bire \cite{Bire01} evaluate Indonesia's Nusa Tenggara Timus province using Malmquist-DEA considering three inputs (namely, number of accommodations, number of restaurants, and number of attractions) and an output (tourist visits) rethinking a new scenario for the regional tourism stakeholders. Pérez León et al. \cite{Perez01} propose an index for measuring tourist destinations in the Caribbean Region, considering 27 indicators in 4 sub-indexes using DEA and goal programming to build composite indicators and measure the competitiveness of destinations. Flegl et al. \cite{Flegl01} measure the hospitality in Mexico using the CCR-DEA model and an input (number of rooms per hotel's star) and three outputs (occupancy rate, tourists arrivals, and related revenue per available room) getting high-efficiency results for national tourism and low-efficiency for international tourism and highlighting that the first is located in land-states and the second in coastal states.

%%%%%%%%%%%%%%%-------------------------------------------------------------------
\begin{table}[t]
\caption{Description and data of the Input and Outputs for the Tourism application}
\label{TableEx0}
%\centering
%\hspace*{-2. cm}
\renewcommand*{\arraystretch}{1.1}
\scalebox{0.73}{
\begin{tabular}{ llp{2.25cm}l }
\toprule
Dimension & Variable & Input/Output & Source
\\
  \hline
\multirow{3}{*}{Economic} & Bed-places (BP) & Input$^\ast$ & Eurostat \& Regional Data Sources  \\
 & Receipts (RCP) & Output &  National Data Sources  \\
 & Overnights (ON) & Output &  National Data Sources  \\ \cline{2-4}
\multirow{2}{*}{Social} & No. Tourism Male Employees (ME) & Output & Eurostat  \\
 & No. Tourism Female Employees (FE) & Output & Eurostat \\\cline{2-4}
Ambiental & Greenhouse Gases (GHG) & Undesirable Output & Eurostat  \\
  \bottomrule
\multicolumn{4}{l}{${}^\ast$ \textit{\footnotesize Interval input obtained from two different databases}}\\%[2ex]
\end{tabular}
}
\renewcommand*{\arraystretch}{1.}
\scalebox{0.75}{
\begin{tabular}{r  c ccccc}
        \toprule
       \multirow{2}{*}{DMUs - Region} & Input & \multicolumn{5}{c}{Outputs} \\
        \cmidrule(lr){2-2} \cmidrule(lr){3-7}
         & BP$^*$ & RCP & ON & FE & ME & GHG \\
          \hline
          Attiki & $[ 62.9 \ , \ 77.41 ]$ & $2\,591.8$ & $4\,973.99$ & $53.0$ & $61.3$ & $19.27$ \\
          Nisia Aigaiou, Kriti & $[ 187.6 \ , \ 242.71 ]$  & $3\, 600.9$ & $7\, 765.65$ & $18.5$ & $21.4$ & $85.85$ \\
          Cataluña & $[ 607.78 \ , \ 791.73 ]$ & $21\, 318.8$ & $20\, 717.24$ & $249.8$ & $215.8$ & $387.72$ \\
          Comunitat Valenciana & $[ 399.66 \ , \ 393.11 ]$ & $9\, 553.1$ & $15\, 830.79$ & $146.5$ & $126.6$ & $221.67$ \\
          Illes Balears & $[ 467.73 \ , \ 443.02 ]$ & $14\,843.4$ & $8\,439.95$ & $78.5$ & $67.9$ & $222.78$ \\
         Provence-Alpes-Côte d'Azur 
         & $[ 677.73 \ , \ 616.56 ] $ & $11\, 779.4$ & $17\, 113.03$ & $68.3$ & $74.7$ & $256.26$ \\
          Jadranska Hrvatska & $[ 1\,080.96 \ , \ 1\, 170.85 ]$ & $2\, 808.9$ & $12\, 219.19$ & $45.4$ & $34.9$ & $44.34$ \\
          Veneto & $[794.25 \ , \ 794.25 ] $ & $29\,396.5$ & $6\, 858.77$ & $257.9$ & $250.9$ & $504.59$ \\
          Campania & $[ 225.17 \ , \ 225.17 ] $ & $4\, 662.9$ & $4\, 878.87$ & $190.5$ & $228.5$ & $156.91$ \\
          Sicilia & $[ 210.92 \ , \ 210.92 ] $ & $3\, 294.4$ & $5\, 802.64$ & $149.7$ & $179.6$ & $127.94$ \\
          Cyprus & $[ 90.19 \ , \ 90.19 ] $ & $3\, 172.1$ & $4\, 241.27$ & $19.6$ & $18.3$ & $90.5$ \\
          Malta & $[ 48.10 \ , \ 52.67 ] $ & $2\, 149.4$ & $3\, 212.46$ & $6.8$ & $12.7$ & $65.08$ \\
        \hline
        Unit of measurement & $10^3$ & $10^6$ \euro & $10^3$ & $10^3$ & $10^3$ & $10^3$ Tonne \\
        \bottomrule
\end{tabular}}
%%%%%%%%%%
\end{table}
%%%%%%%%%%%%%%%-------------------------------------------------------------------

% %%%%%%%%%%%%%%%%%
\begin{figure}[t]
    \centering
    %\hspace*{-3.2 cm}
    \includegraphics[width=0.95\textwidth]{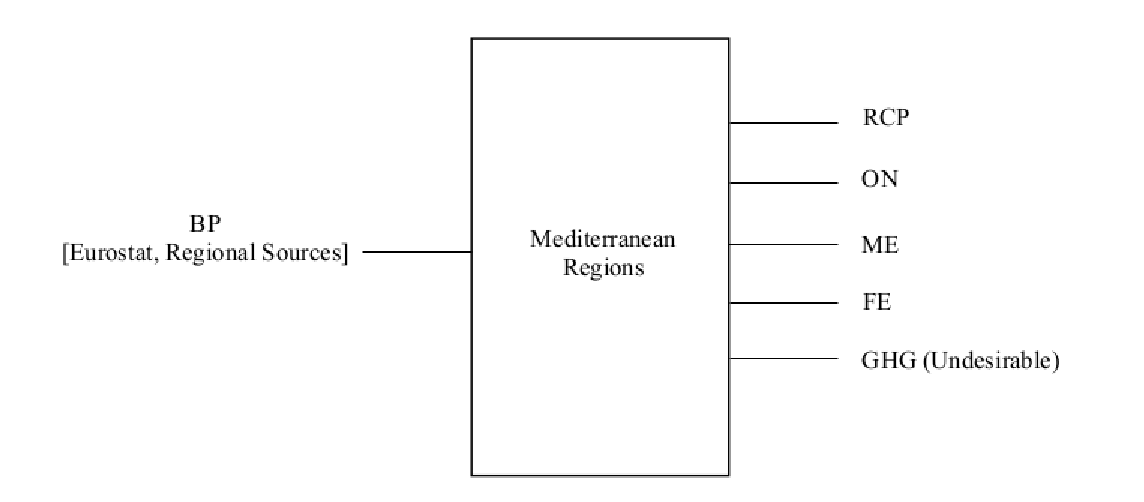}
    \caption{Input and Outputs considered in this sustainable tourism application. }
    \label{fig:Inputs_outputs_box}
\end{figure}

 %%%%%%%%%%%%%%%

\subsection{Variables and Data}

 Several input and output variables are considered from the three sustainable dimensions. See Table \ref{TableEx0} and Figure \ref{fig:Inputs_outputs_box}. The data refer to the year 2019, prior to the pandemic. We have used the Eurostat database and Regional databases from different Mediterranean regions. We would have liked to include additional regions and variables, but this was prevented by data availability. As a novelty in tourism studies, the variable Bed-places is considered an interval, estimated as a confidence interval from the available data. Thus, Bed-places data come from two databases with different data for some regions. Also, similar to other studies like \cite{Guo02} and \cite{Gao01}, GHG emissions have been considered undesirable. Using interval variables, the proposed DEA model is a novelty in sustainability tourism efficiency assessment.

%%%%%%%%%%%%%%%-------------------------------------------------------------------
\begin{sidewaystable}[p]
%\begin{table}
    \centering
    \caption{Results from the (EIMIL) \eqref{prob: EIMIL} (second column) and (S-EIMIL) \eqref{prob: SEIMIL} (third column) models, respectively. The corresponding ranking is given in the fourth column. We also include the slacks and the targets. Only the inefficient DMUs have non-zero slacks. For clarity, we represent those null interval slacks as zero, i.e., $0 \equiv [0,0]$.
    }
\label{TableEx2.a}
    \hspace*{-0.5cm}
    \scalebox{0.7}{
    \begin{tabular}{|r| c |c | c | cc|cc| cc|cc|cc|cc|}
    \hline
       \multirow{2}{*}{DMU} &\multirow{2}{*}{$EI(X_p,Y_p)$}& 
       \multirow{2}{*}{$SEI(X_p,Y_p)$}&\multirow{2}{*}{Ranking}&
       \multicolumn{2}{c|}{Input slacks} & \multicolumn{10}{c|}{Output slacks} \\
        \cline{5-16}
         &&&&$sl^x_1$ & $su^x_1$ &
         $sl^y_1$ & $su^y_1$ & $sl^y_2$ & $su^y_2$ &
         $sl^y_3$ & $su^y_3$ & $sl^y_4$ & $su^y_4$ 
         & $sl^y_5$ & $su^y_5$\\
          \hline
         Attiki  &  0  &  4.042  &  1  & $0$ & $0$ & $0$ & $0$ & $0$ & $0$ & $0$ & $0$ & $0$ & $0$ & $0$ & $0$ \\
        Nisia Aigaiou, Kriti  &  7.897  &  $--$  &  12  & $[ 22.57 , 66.74 ]$ & $0$ & $0$ & $1\,003.21$ & $0$ & $0$ & $72.97 $ & $0$ & $74.19$ & $0$ & $0$ & $0$ \\
        Cataluña  &  0  &  0.699  &  2  & $0$ & $0$ & $0$ & $0$ & $0$ & $0$ & $0$ & $0$ & $0$ & $0$ & $0$ & $0$ \\
        Comunitat Valenciana  &  0  &  0.228  &  7  & $0$ & $0$ & $0$ & $0$ & $0$ & $0$ & $0$ & $0$ & $0$ & $0$ & $0$ & $0$ \\
        Illes Balears  &  0  &  0.089  &  8  & $0$ & $0$ & $0$ & $0$ & $0$ & $0$ & $0$ & $0$ & $0$ & $0$ & $0$ & $0$ \\
        Provence-Alpes-Côte d'Azur  &  2.556  &  $--$  &  10  & $0$ & $[ 0.62 , 27.83 ]$ & $0$ & $1\,182.26 $ & $0$ & $0$ & $ 101.23 $ & $0$ & $71.07 $ & $0$ & $0$ & $0$ \\
        Jadranska Hrvatska  &  0  &  0.483  &  3  & $0$ & $0$ & $0$ & $0$ & $0$ & $0$ & $0$ & $0$ & $0$ & $0$ & $0$ & $0$ \\
        Veneto  &  0  &  0.446  &  5  & $0$ & $0$ & $0$ & $0$ & $0$ & $0$ & $0$ & $0$ & $0$ & $0$ & $0$ & $0$ \\
        Campania  &  0  &  0.455  &  4  & $0$ & $0$ & $0$ & $0$ & $0$ & $0$ & $0$ & $0$ & $0$ & $0$ & $0$ & $0$ \\
        Sicilia  &  0  &  0.003  &  9  & $0$ & $0$ & $0$ & $0$ & $0$ & $0$ & $0$ & $0$ & $0$ & $0$ & $0$ & $0$ \\
        Cyprus  &  4.221  &  $--$  &  11  & $0$ & $[ 0 , 12.04 ]$ & $0$ & $0$ & $0$ & $406.38 $ & $0$ & $ 28.75 $ & $37.46 $ & $0$ & $ 49.3 $ & $0$ \\
        Malta  &  0  &  0.392  &  6  & $0$ & $0$ & $0$ & $0$ & $0$ & $0$ & $0$ & $0$ & $0$ & $0$ & $0$ & $0$ \\
         \hline
         %%%%%%%%%%%%%%%%%%%%%%%%%%%%%%%%%%%%%%%%%%%%%%%%%%%%%%%
        \multicolumn{4}{c|}{}&
        \multicolumn{2}{c|}{Input Target } & \multicolumn{10}{c|}{Output Targets} 
        \\ \cline{5-16}
        \multicolumn{4}{r|}{DMU}&
         \multicolumn{2}{c|}{$x^{target}_{1p}$ } &
         \multicolumn{2}{c|}{$y^{target}_{1p}$ }&
         \multicolumn{2}{c|}{$y^{target}_{2p}$ }&
         \multicolumn{2}{c|}{$y^{target}_{3p}$ } &
         \multicolumn{2}{c|}{$y^{target}_{4p}$ } &
        \multicolumn{2}{c|}{$y^{target}_{5p}$ } 
        \\ \cline{2-16}
        \multicolumn{4}{r |}{  Attiki } & \multicolumn{2}{c|}{$[ 62.9 \ , \  77.41 ]$} & \multicolumn{2}{c|}{$2\,591.8 $} &\multicolumn{2}{c|}{$4\,973.99 $} & \multicolumn{2}{c|}{$53$} & \multicolumn{2}{c|}{$61.3$} & \multicolumn{2}{c|}{$19.27 $} \\
        \multicolumn{4}{r |}{  Nisia Aigaiou, Kriti } & \multicolumn{2}{c|}{$[ 165.03 \ , \ 175.97 ]$} & \multicolumn{2}{c|}{$ 4\,604.11 $} &\multicolumn{2}{c|}{$7\,765.65 $} & \multicolumn{2}{c|}{$91.47 $} & \multicolumn{2}{c|}{$95.59 $} & \multicolumn{2}{c|}{$5.85$} \\
        \multicolumn{4}{r |}{  Cataluña } & \multicolumn{2}{c|}{$[ 607.78 \ , \ 791.73 ]$} & \multicolumn{2}{c|}{$ 21\,318.8 $} &\multicolumn{2}{c|}{$ 20\,717.24 $} & \multicolumn{2}{c|}{$ 249.8 $} & \multicolumn{2}{c|}{$ 215.8 $} & \multicolumn{2}{c|}{$ 387.72 $} \\
        \multicolumn{4}{r |}{  Comunitat Valenciana } & \multicolumn{2}{c|}{$[ 393.11 \ , \ 399.66 ]$} & \multicolumn{2}{c|}{$9\,553.1 $} &\multicolumn{2}{c|}{$15\,830.79 $} & \multicolumn{2}{c|}{$146.5 $} & \multicolumn{2}{c|}{$126.6 $} & \multicolumn{2}{c|}{$ 221.67 $} \\
        \multicolumn{4}{r |}{  Illes Balears } & \multicolumn{2}{c|}{$[ 443.02 \ , \ 467.73 ]$} & \multicolumn{2}{c|}{$  14\,843.4 $} &\multicolumn{2}{c|}{$ 8\,439.94 $} & \multicolumn{2}{c|}{$ 78.5 $} & \multicolumn{2}{c|}{$ 67.9 $} & \multicolumn{2}{c|}{$ 222.78 $} \\
        \multicolumn{4}{r |}{  Provence-Alpes-Côte d'Azur } & \multicolumn{2}{c|}{$[ 588.72 \ , \ 677.12 ]$} & \multicolumn{2}{c|}{$ 12\,961.66 $} &\multicolumn{2}{c|}{$ 17\,113.03 $} & \multicolumn{2}{c|}{$ 169.53 $} & \multicolumn{2}{c|}{$ 145.77 $} & \multicolumn{2}{c|}{$ 256.26 $} \\
        \multicolumn{4}{r |}{  Jadranska Hrvatska } & \multicolumn{2}{c|}{$[ 1\,080.96 \ , \ 1\,170.85 ]$} & \multicolumn{2}{c|}{$ 2\,808.9 $} &\multicolumn{2}{c|}{$ 12\,219.19 $} & \multicolumn{2}{c|}{$ 45.4 $} & \multicolumn{2}{c|}{$ 34.9 $} & \multicolumn{2}{c|}{$ 44.34 $} \\
        \multicolumn{4}{r |}{  Veneto } & \multicolumn{2}{c|}{$[ 794.25 \ , \ 794.25 ]$} & \multicolumn{2}{c|}{$ 29\,396.5 $} &\multicolumn{2}{c|}{$ 6\,858.77 $} & \multicolumn{2}{c|}{$ 257.9 $} & \multicolumn{2}{c|}{$ 250.9 $} & \multicolumn{2}{c|}{$ 504.59 $} \\
        \multicolumn{4}{r |}{  Campania } & \multicolumn{2}{c|}{$[ 225.17 \ , \ 225.17 ]$} & \multicolumn{2}{c|}{$ 4\,662.9 $} &\multicolumn{2}{c|}{$ 4\,878.87 $} & \multicolumn{2}{c|}{$ 190.5 $} & \multicolumn{2}{c|}{$ 228.5 $} & \multicolumn{2}{c|}{$ 156.91 $} \\
        \multicolumn{4}{r |}{  Sicilia } & \multicolumn{2}{c|}{$[ 210.92 \ , \ 210.92 ]$} & \multicolumn{2}{c|}{$ 3\,294.4 $} &\multicolumn{2}{c|}{$ 5\,802.64 $} & \multicolumn{2}{c|}{$ 149.7 $} & \multicolumn{2}{c|}{$179.6 $} & \multicolumn{2}{c|}{$ 127.94 $} \\
        \multicolumn{4}{r |}{  Cyprus } & \multicolumn{2}{c|}{$[ 78.15 \ , \ 90.19 ]$} & \multicolumn{2}{c|}{$ 3\,172.1 $} &\multicolumn{2}{c|}{$ 4\,647.65 $} & \multicolumn{2}{c|}{$ 48.35 $} & \multicolumn{2}{c|}{$ 55.76 $} & \multicolumn{2}{c|}{$ 41.2 $} \\
        \multicolumn{4}{r |}{  Malta } & \multicolumn{2}{c|}{$[48.1 \ , \ 52.67 ]$} & \multicolumn{2}{c|}{$ 2\,149.4 $} &\multicolumn{2}{c|}{$ 3\,212.45 $} & \multicolumn{2}{c|}{$ 6.8 $} & \multicolumn{2}{c|}{$ 12.7 $} & \multicolumn{2}{c|}{$ 65.08 $} \\
             \cline{2-16}
    \end{tabular}
  }  
\end{sidewaystable}
%%%%%%%%%%%%

%%%-------------------------------------------------------------------
\begin{table}
    \centering
    \caption{Comparison with other models from the Literature. }
\label{TableEx2.a}
    %\hspace*{-2.0cm}
    %\vspace*{-1.cm}
    \scalebox{0.8}{
    \begin{tabular}{|r| c c c |c c | }
    \hline
       \multirow{2}{*}{DMU} &
        \multicolumn{3}{c|}{This work} &
       \multicolumn{2}{c|}{Azizi et al. (2015)}
       \\[1ex]
         &  $EI(X_p,Y_p)$ &  $SEI(X_p,Y_p)$ & $ Ranking $ &
         $[\Phi_p^L,\Phi_p^U]$ & $ Ranking $ \\
        \hline
         Attiki  &  0  &  4.042  &  1  &
         $[ 2.513 , 2.889 ]$ & 1
         \\
        Nisia Aigaiou, Kriti  &  7.897  &  $--$  &  12  &  
        $[ 0.406 , 0.469 ]$ & 7
        \\
        Cataluña  &  0  &  0.699  & 2  & 
         $[ 0.736 , 0.882 ]$ & 9
        \\
        Comunitat Valenciana  &  0  & 0.228  &  7  &
         $[ 0.874 , 0.884 ]$ & 4
        \\
        Illes Balears  &  0  & 0.089  &  8  &  
         $[ 0.434 , 0.448 ]$  & 3
       \\
        Provence-Alpes-Côte d'Azur  &  2.556  &  $--$  &  10  &  
        $[ 0.448 , 0.472 ]$ & 12
       \\
        Jadranska Hrvatska  &  0  &  0.483  &  3  &  
         $[ 0.981 , 1.04]$ & 10
       \\
        Veneto  &  0  & 0.446  &  5  & 
         $[ 0.431 , 0.431 ]$ & 11
        \\
        Campania  &  0  &  0.455  &  4  & 
        $[ 1, 1]$  & 6
       \\
        Sicilia  &  0  & 0.003  &  9  & 
        $[ 0.619 , 0.619 ]$ & 5
        \\
        Cyprus  &  4.221  & $--$  &  11  & 
         $[ 0.476 , 0.476 ]$ & 2
        \\
        Malta  &  0  &  0.392  &  6  &  
         $[ 0.44, 1.023 ]$ & 8
        \\
         \hline
     \end{tabular}}
    \end{table}

%%%-------------------------------------------------------------------
\subsection{Results and discussion}

The inefficiency scores of the proposed (EIMIL) DEA model $EI(X_p,Y_p)$ with its corresponding targets and slacks intervals are shown in Table \ref{TableEx2.a}. 
Because of the relatively small dataset, only three DMUs are inefficient: Nisia Aigaiou-Kriti, Cyprus, and the Provence-Alpes-Côte d’Azur. That is why, in order to rank the DMUs fully, the super SBI interval DEA model (S-PEIMIL) has been solved. The corresponding super-inefficiency scores $SEI(X_p,Y_p)$ and the final ranking are also shown. 

Note that, regarding sustainable tourism, the two Greek regions considered in this studio have disparate results, with Attiki and Nisia Aigaiou (Kriti) at the top and bottom of the ranking, respectively. The second-best score is for Cataluña, with similar scores as the Croatian region Jadranska Hrvatska, followed by the Italian regions of Campania and Veneto. On the bottom side, Cyprus and the Provence-Alpes-Côte d’Azur, in addition to Kriti, are inefficient. 

Regarding inefficient regions, Nisia Aigaiou (Kriti) needs to increase receipts by around 28\% and female and male employment by up to 494\% and 447\%, respectively, to reach the frontier of tourism sustainability. Regarding the number of beds, the interval variable has an excess of 28\% and 12\% with respect to their current value according to the Eurostat and the regional databases, respectively. Similarly, Provence-Alpes-Côte d'Azur has margins for improvement of 248\% and 195\% for female and male employment, respectively, and for increasing the receipts by 10\%. The interval input variable shows zero slacks for the Eurostat database and around 5\% slack for the regional database. Finally, Cyprus is the only study area with a margin of improvement in greenhouse gas emissions. Namely, they could decrease by 46\%. It also has margins of improvement of 246\% and 305\% in female and male employment, respectively, and of 10\% in overnight stays. As for the BP interval input, the slack is zero for the Eurostat database and around 5\% for the regional database.  

For comparison purposes, Table \ref{TableEx2.a} shows the overall efficiency intervals $[\Phi_p^L,\Phi_p^U]$ computed by the non-oriented SBM  model for interval data proposed by Azizi et al. \cite{azizi}. Their method computes an interval measure of efficiency although it does not compute targets. Their approach uses a preference-degree approach for comparing and ranking the DMUs. The corresponding ranking computed by these authors for this dataset is also included in the Table. Their ranking and the proposed approach are not correlated (Spearman correlation coefficient=-0.056). This is undoubtedly due to their considering a double frontier approach.

\begin{figure*}
    \centering
    \vspace*{-2.cm}
    \includegraphics[width = .75\textwidth]{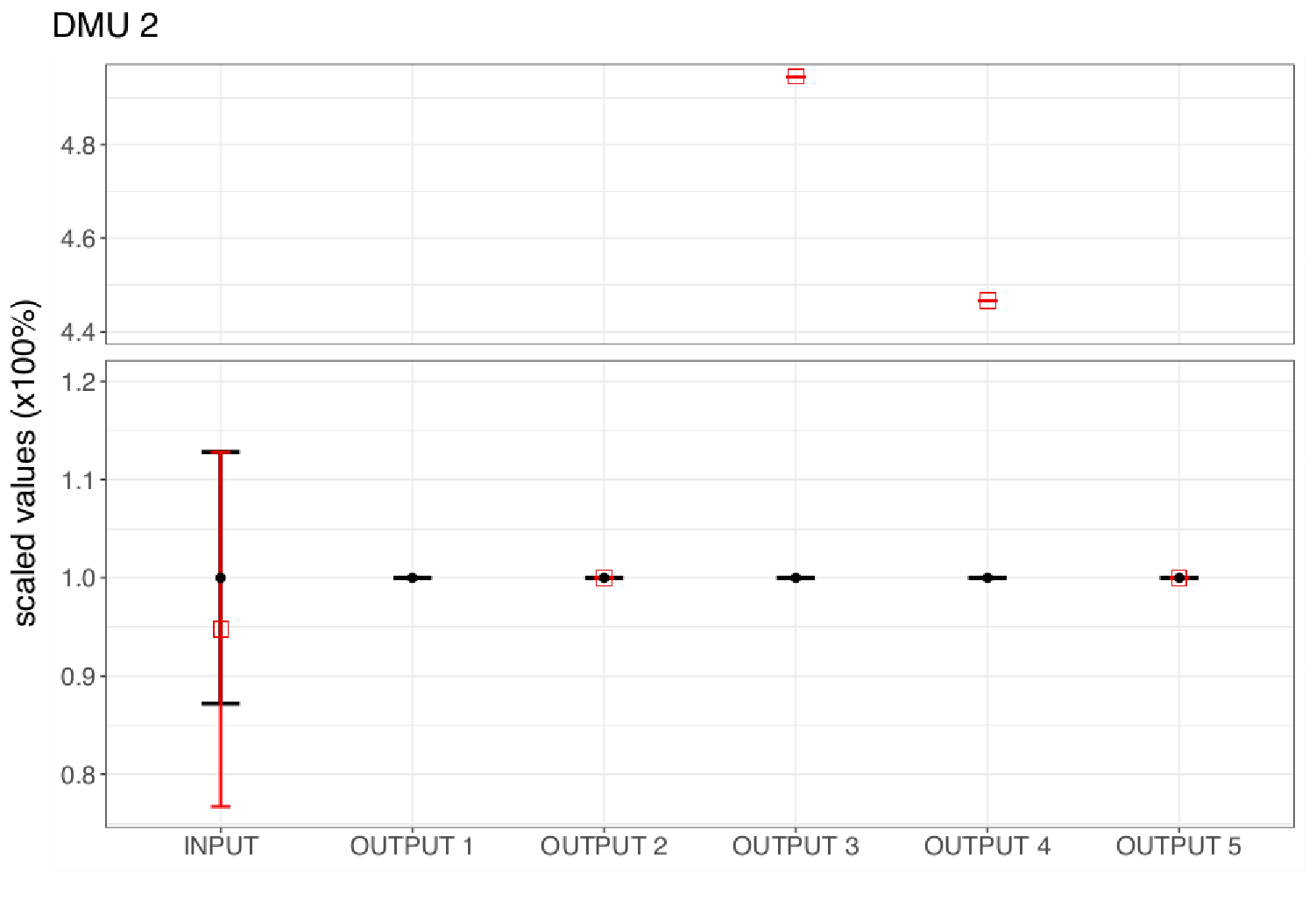}
    \includegraphics[width = .75\textwidth]{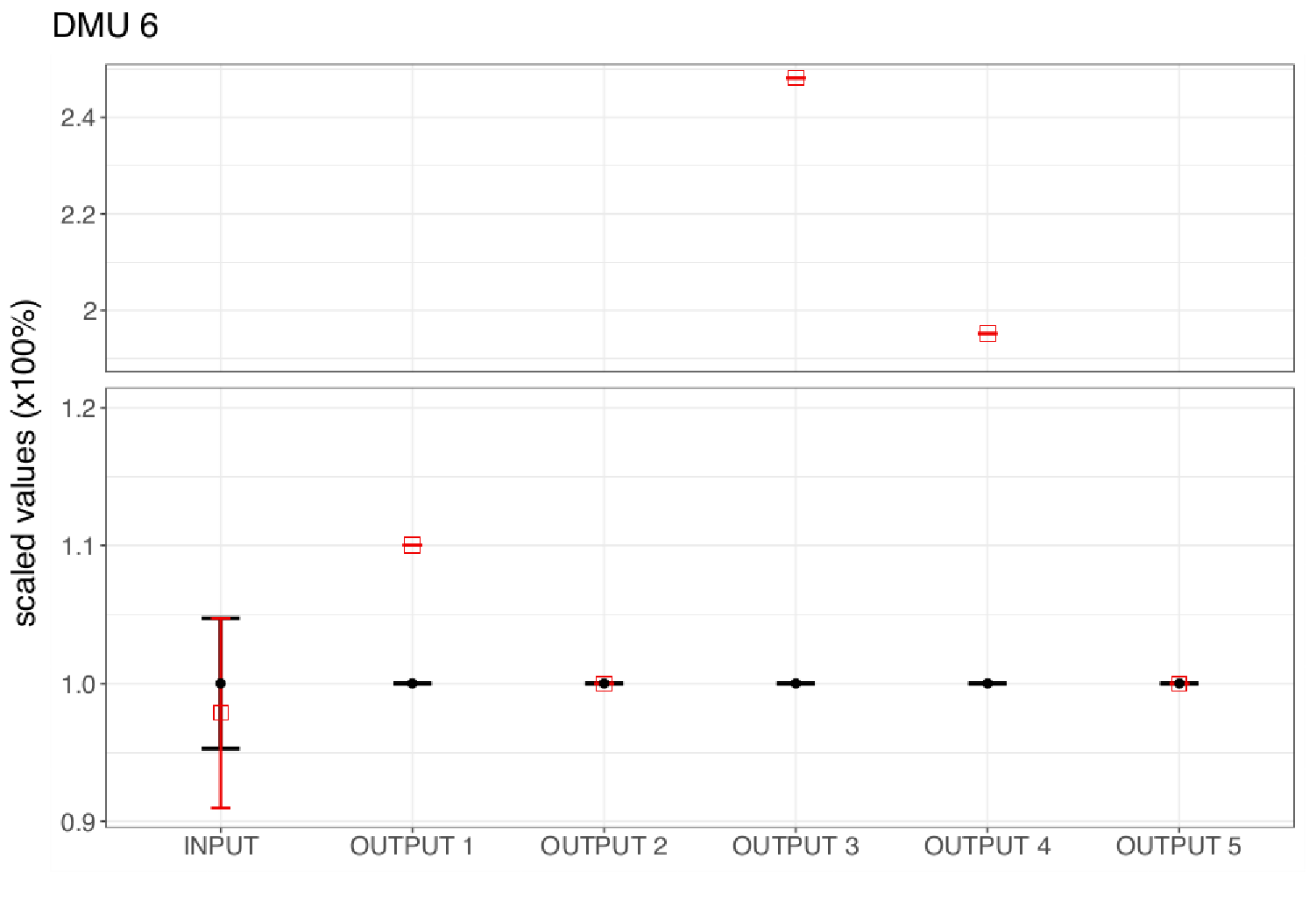}
    \includegraphics[width = .75\textwidth]{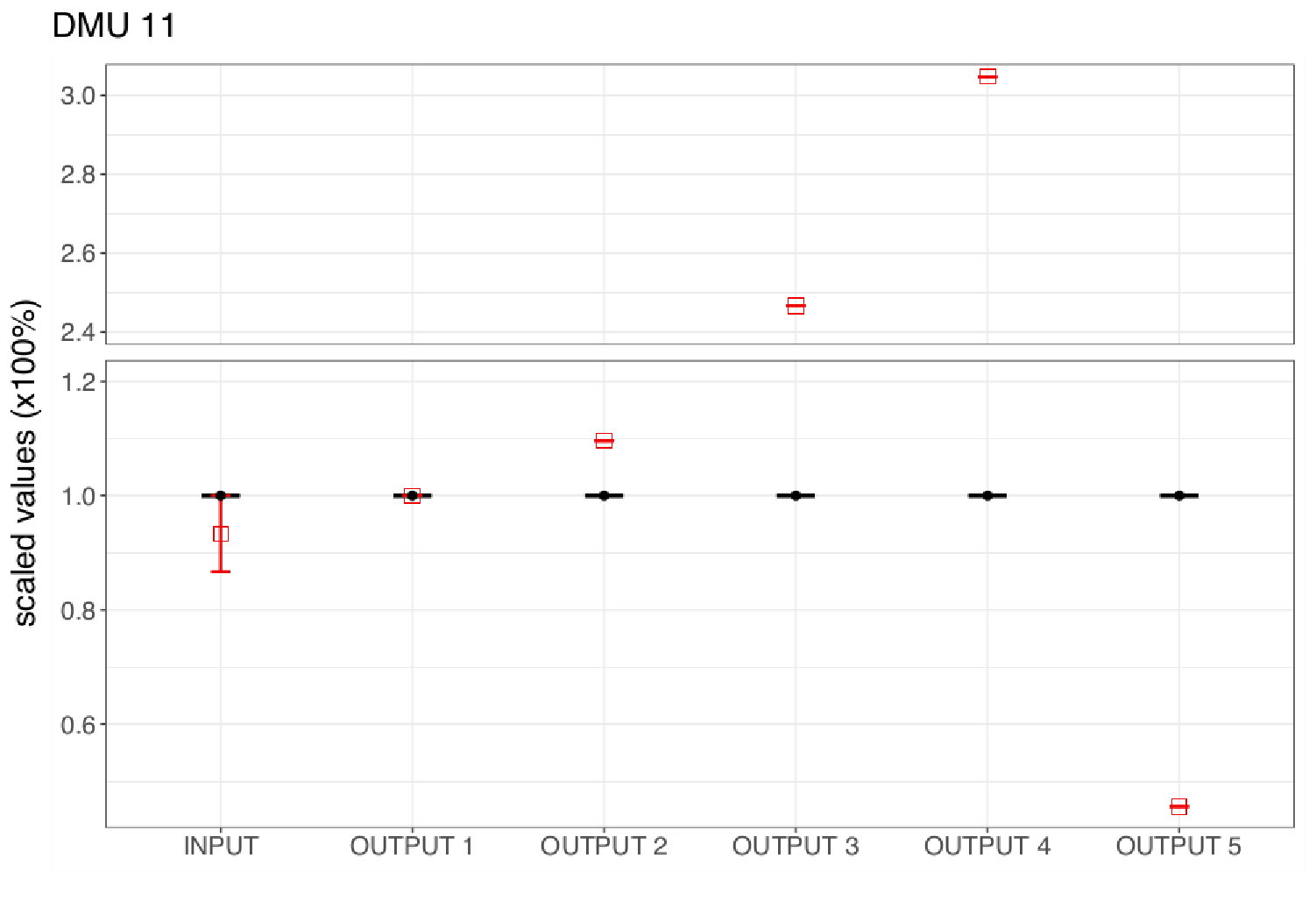}
        \caption{Input and outputs targets for the three inefficient DMUs: Nisia Aigaiou, Kriti  (labeled DMU 2), Provence-Alpes-Côte d'Azur (labeled DMU 6), and Cyprus (labeled DMU 11).}
        \label{fig:ineff dmus}
\end{figure*}

Figure \ref{fig:ineff dmus} shows the observed and target input and output targets for these three inefficient DMUs. The values are scaled by the corresponding observed data to facilitate their
comparison. Note that, as the fifth output is undesirable, the corresponding targets stayed the same or were reduced. Note also that, in this application, only the input variable involves interval data.

%%%%%%%%%%%%%%%-------------------------------------------------------------------
%%%%%%%%%%%%%%%-------------------------------------------------------------------

\section{Conclusions}\label{conclusions}
%%%%%%%%%%%%%

This paper has proposed a new interval-valued DEA approach and associated slacks-based inefficiency measures.  It requires solving a mixed-integer linear program that allows computing the corresponding input and output targets. A super-efficiency version of the model has also been formulated in case fully ranking the DMUs is desired.

An application for sustainable tourism efficiency assessment has been presented. The input and output variables span the three sustainability dimensions, including the environmental dimension, represented by GHG emissions from tourism activities. The need to apply an interval DEA approach comes from the fact that for the input variable (Bed-places), the data comes from two different data sources, and in some cases, the corresponding values do not coincide. This is something that often occurs in practice. In order to avoid loss of information, it was decided to represent that variable as an interval using the values from the two sources as limits. The proposed approach has handled this type of variable computing inefficiency scores and targets for all the DMUs. In this application, given the small dataset available, many DMUs were labeled as efficient, which also required solving the corresponding super-efficiency DEA model.

As a continuation of this research, we will apply this approach to other sectors (e.g., healthcare, sports, etc.) where input and output interval data can occur. Theoretically, we could extend it to Network DEA scenarios, i.e., production systems involving multiple interconnected processes.

\section{{Acknowledgements}} 

The first and second authors are partially supported by grant PID2019-105824GB-I00. The fifth author acknowledges the financial support
of the Spanish Ministry of Science and Innovation, grant PID2021-124981NB-I00.

\appendix

%%%%%%%%%%%%%%%%%%%%%%%
\end{document}